\documentclass[12pt,english]{amsart}
\usepackage[charter]{mathdesign}

\usepackage[T1]{fontenc}
\usepackage[latin9]{inputenc}
\usepackage[a4paper]{geometry}
\usepackage[active]{srcltx}
\usepackage{verbatim}
\usepackage{textcomp}
\usepackage{mathtools}
\usepackage{amsthm}
\usepackage{stmaryrd}
\usepackage{setspace}
\makeatletter
\theoremstyle{plain}
\newtheorem{thm}{\protect\theoremname}[section]
  \theoremstyle{plain}
  \newtheorem{cor}[thm]{\protect\corollaryname}
  \theoremstyle{definition}
  \newtheorem{defn}[thm]{\protect\definitionname}
  \theoremstyle{remark}
  \newtheorem{rem}[thm]{\protect\remarkname}
  \theoremstyle{plain}
  \newtheorem{prop}[thm]{\protect\propositionname}
  \theoremstyle{definition}
  \newtheorem{example}[thm]{\protect\examplename}

\def\sym{\mathsf{Sym}}

\def\P{\mathsf{PE}}

\usepackage{babel}

\title[Plethystic calculus and permutations]
{Plethystic exponential calculus \\ and characteristic polynomials of permutations}

\author[C. Florentino]{Carlos Florentino}

\address{Departamento de Matem\'{a}tica, Faculdade de Ci\^{e}ncias, Univ. de Lisboa, 
Edf. C6, Campo Grande 1749-016 Lisboa, Portugal}

\email{caflorentino@fc.ul.pt}

\thanks{Partially supported from the
Project QuantumG PTDC/MAT-PUR/30234/2017, FCT, Portugal}

\@namedef{subjclassname@2020}{\textup{2020} Mathematics Subject Classification}


\usepackage{babel}
\usepackage[all]{xy}

\renewcommand{\textendash}{--}

\makeatother

\usepackage{babel}
  \providecommand{\corollaryname}{Corollary}
  \providecommand{\definitionname}{Definition}
  \providecommand{\examplename}{Example}
  \providecommand{\propositionname}{Proposition}
  \providecommand{\remarkname}{Remark}
\providecommand{\theoremname}{Theorem}

\begin{document}
\begin{abstract}
We prove a family of identities, expressing generating functions of
powers of characteristic polynomials of permutations, as finite or
infinite products. %
These generalize formulae first obtained in a study of the geometry/topology
of symmetric products of real/algebraic tori. The proof uses formal
power series expansions of plethystic exponentials, and has been motivated
by some recent applications of these combinatorial tools in supersymmetric
gauge and string theories. Since the methods are elementary, we tried
to be self-contained, and relate to other topics such as the $q$-binoomial
theorem, and the cycle index and Molien series for the symmetric group.
\end{abstract}

\maketitle

\section{Introduction}

\thispagestyle{empty}

The study of infinite products is a classical theme, and such objects
arise in many parts of mathematics since the early works of L. Euler.
They are fundamental in the theory of holomorphic functions of a complex
variable, and some notable infinite products of simple fractions are
generating functions for counting various types of combinatorial objects,
such as partitions of natural integers. 

Power series expansions of infinite products are also ubiquitous,
and some have become famous, such as Euler's pentagonal number theorem,
or Jacobi's triple product identity; these have a wide range of applications
in arithmetic geometry, complex function theory, and number theory,
\emph{et cetera}. However, it is not so often that the power series
expansion of a given infinite product can be written only in terms
of elementary functions and concepts. 

The purpose of this article is two-fold. Firstly, we want to provide
a completely elementary proof of the following family of simple identities,
parametrized by an integer $r\in\mathbb{Z}$. Let $\mathbb{Z}\llbracket x,y\rrbracket$
be the ring of formal power series in the variables $x$ and $y$
with integer coefficients. We denote by $S_{n}$ the symmetric group
on $n$ letters, $n\in\mathbb{N}$. 
\begin{thm}
\label{thm:Main} Fix $r\in\mathbb{Z}$. Then, the following are equalities
in $\mathbb{Z}\llbracket x,y\rrbracket$
\begin{equation}
1+\sum_{n\geq1}\sum_{\sigma\in S_{n}}\frac{y^{n}}{n!}\det(I_{n}-xM_{\sigma})^{r}=\ \prod_{k\geq0}(1-yx^{k})^{-\binom{k-r-1}{k}},\label{eq:main-identity}
\end{equation}
where $M_{\sigma}$ is the $n\times n$ permutation matrix of $\sigma\in S_{n}$,
and $I_{n}$ the identity matrix of the same size.
\end{thm}

The above identities are actually valid in $\mathbb{R}\llbracket x,y\rrbracket$
(and also in $\mathbb{C}\llbracket x,y\rrbracket$, with appropriate
use) upon letting $r$ be an arbitrary real number, with the binomial
numbers defined by $\binom{r}{k}=\frac{r(r-1)\cdots(r-k+1)}{k!}$,
for all $k\geq0$ (and $\binom{0}{0}=1$). In the important cases
when $r$ is a \emph{non-negative integer}, the right hand side product
becomes finite:
\begin{equation}
1+\sum_{n\geq1}\sum_{\sigma\in S_{n}}\frac{y^{n}}{n!}\det(I_{n}-xM_{\sigma})^{r}=\prod_{k=0}^{r}(1-yx^{k})^{(-1)^{k+1}\binom{r}{k}},\label{eq:positive-r}
\end{equation}
given that $\binom{k-r-1}{k}=(-1)^{k}\binom{r}{k}$, for $r\in\mathbb{N}$;
and the series converges absolutely for $x,y\in\mathbb{C}$, with
$|x|<1$, $|y|<2^{-r}$ (and uniformly in compact subsets).

For a given $\sigma\in S_{n}$, the polynomials $q_{\sigma}(x):=\det(I_{n}-xM_{\sigma})$
will be called \emph{characteristic polynomials of permutations} (the
more common polynomials $\pm\det(I_{n}\lambda-M_{\sigma})$ are obtained
by an obvious substitution). Hence, the left hand side of the identities
above are (ordinary)\emph{ generating functions for the average $r$th
powers of} \emph{characteristic polynomials of permutations}.

The proof of Theorem \ref{thm:Main} is presented in Section 3 using
only a few simple and well-known interesting tools in combinatorics,
borrowed from the theory of symmetric functions: the so-called plethystic
exponential calculus. Hence, our second purpose is to provide a concise
introduction, in Section 2, to the aspects of plethystic calculus
showing up in the proof - the plethystic exponential in a ring of
formal power series - to help making them better known. We believe
that this exponential is bound to play an increasing role in many
topics, as it is shown by several recent and diverse applications
in so many areas of Mathematics and Mathematical-Physics (see, for
example, \cite{BBS,BFHH,DMVV,FHH}).

Actually, to obtain our main identities, we do not need the full power
of plethysms from the theory of symmetric functions or lambda rings,
and we only need to establish product-sum identities for certain plethystic
exponentials. For $r\geq1$, denote by $\P(f(\underline{x},y))$ the
plethystic exponential of a formal power series $f$, with $\mathbb{R}$
coefficients, in $r+1$ variables $\underline{x}=(x_{1},\cdots,x_{r})$
and $y$, without constant term (see Definition \ref{def:PE-2var}
below). Let us use the multi-index notation for the $\underline{x}$
variables: $\mathbf{k}=(k_{1},\cdots,k_{r})\in\mathbb{N}_{0}^{r}$,
($\mathbb{N}_{0}=\mathbb{N}\cup\{0\}$) and write $\mathbf{x}^{\mathbf{k}}=x_{1}^{k}\cdots x_{r}^{k}$. 

Let $g(\underline{x})=\sum_{\mathbf{k}\in\mathbb{N}_{0}^{r}}a_{\mathbf{k}}\mathbf{x}^{\mathbf{k}}\in\mathbb{R}\llbracket\underline{x}\rrbracket$,
so that $a_{\mathbf{k}}=a_{k_{1},\cdots,k_{r}}\in\mathbb{R}$. The
source of our product-sum identities, is the following evaluation
of plethystic exponentials, both as an infinite product, and as a
sum over \emph{partitions} of the natural number $n$:
\begin{equation}
\P\left(g(\underline{x})\,y\right)=\prod_{\mathbf{k}\in\mathbb{N}_{0}^{r}}(1-y\,\mathbf{x}^{\mathbf{k}})^{-a_{\mathbf{k}}}=1+\sum_{n\geq1}\left(\sum_{\lambda\vdash n}\ \prod_{j=1}^{n}\frac{g(\underline{x}^{j})^{\lambda(j)}}{\lambda(j)!\,j^{\lambda(j)}}\right)\,y^{n}\,\in\mathbb{R}\llbracket\underline{x},y\rrbracket.\label{eq:product-sum}
\end{equation}
Above, we write $\lambda\vdash n$ when $\lambda$ is a partition
of $n$, and $\lambda(j)\geq0$ denotes the number of parts of $\lambda$
of size equal to $j\in\{1,\cdots,n\}$. We also used the notation:
$\underline{x}^{j}:=(x_{1}^{j},\cdots,x_{r}^{j})$, for $j\in\mathbb{N}$.%
{} 

We will apply Equation \eqref{eq:product-sum} to several polynomials/formal
power series $g(\underline{x})$. By using appropriate rational functions
in $r+s$ variables, now denoted $x_{1},\cdots,x_{r}$ and $q_{1},\cdots,q_{s}$
(given the interesting relations to $q$-series, see §3.2) we obtain
the following, returning now to $\mathbb{Z}$ coefficients.
\begin{thm}
\label{thm:multivar} Fix $r,s\in\mathbb{N}_{0}=\mathbb{N}\cup\{0\}$.
With the same notations as in Theorem \ref{thm:Main}, we have:
\begin{eqnarray}
1+\sum_{n\geq1}\frac{y^{n}}{n!}\sum_{\sigma\in S_{n}}\prod_{j=1}^{r}\prod_{k=1}^{s}\frac{\det(I_{n}-x_{j}M_{\sigma})}{\det(I_{n}-q_{k}M_{\sigma})} & = & \prod_{\mathbf{j}\in\{0,1\}^{r}}\prod_{\mathbf{k}\in\mathbb{N}_{0}^{s}}(1-y\,\mathbf{x}^{\mathbf{j}}\,\mathbf{q}^{\mathbf{k}})^{(-1)^{|\mathbf{j}|+1}},\label{eq:multivar}
\end{eqnarray}
as identities in $\mathbb{Z}\llbracket x_{1},\cdots,x_{r},q_{1},\cdots,q_{s},y\rrbracket$,
with $|\mathbf{j}|=\sum_{i=1}^{r}j_{i}$. 
\end{thm}

Note that Theorem \ref{thm:Main} follows from Theorem \ref{thm:multivar},
by letting $s$ or $r$ equal zero, and identifying all the remaining
variables. We now present a couple of interesting consequences of
these identities. First note that, from Theorem \eqref{thm:Main},
our basic generating function:
\[
\Phi^{r}(x,y):=1+\sum_{n\geq1}\sum_{\sigma\in S_{n}}\frac{y^{n}}{n!}\det(I_{n}-xM_{\sigma})^{r},
\]
has always a factor $1/(1-y)$, corresponding to a pole of order 1
at $y=1$. Its residue is again expressed as a product, and has the
following enumerative interpretation. %
\begin{cor}
\label{cor:residue}For every $r\in\mathbb{Z}$, we have:
\[
\mathrm{Res}_{y=1}\Phi^{r}(x,y)=-\prod_{k\geq1}(1-x^{k}){}^{-\binom{k-r-1}{k}}.
\]
In particular, when $c=-r>0$, $\mathrm{Res}_{y=1}\Phi^{-c}(x,y)$
is the (the negative of) the generating function for partitions of
a natural number $n$, with parts colored with $c$ distinct colors.
\end{cor}

We also obtain, in Section 3, a recursion formula for the \emph{average}
$r$ power of characteristic polynomials of permutations:
\begin{equation}
Q_{n}^{r}(x):=\frac{1}{n!}\sum_{\sigma\in S_{n}}\det(I_{n}-xM_{\sigma})^{r},\label{eq:char-pol-power}
\end{equation}
and, more generally, with $\mathbf{r}=(r_{1},\cdots,r_{s})\in\mathbb{R}^{s}$,
for any formal power series of the form:
\[
Q_{n}^{\mathbf{r}}(\underline{x}):=\frac{1}{n!}\sum_{\sigma\in S_{n}}\det(I_{n}-x_{1}M_{\sigma})^{r_{1}}\cdots\det(I_{n}-x_{s}M_{\sigma})^{r_{s}}\quad\in\mathbb{R}\llbracket\underline{x}\rrbracket.
\]
\begin{cor}
\label{cor:recursion}Fix $\mathbf{r}=(r_{1},\cdots,r_{s})\in\mathbb{R}^{s}$.
Then, as elements in $\mathbb{R}\llbracket\underline{x}\rrbracket$,
we have:
\[
Q_{n}^{\mathbf{r}}(\underline{x})=\frac{1}{n}\sum_{k=1}^{n}h(\underline{x}^{k})\,Q_{n-k}^{\mathbf{r}}(\underline{x}),
\]
with $h(\underline{x}):=(1-x_{1})^{r_{1}}\cdots(1-x_{s})^{r_{s}}$.
\end{cor}

This recursion relation is an immediate consequence of the interpretation
of such formal power series $Q_{n}^{\mathbf{r}}(\underline{x})$ as
evaluations of the cycle index of the symmetric group $S_{n}$ on
$h(\underline{x})\in\mathbb{R}\llbracket\underline{x}\rrbracket$
(see §3.3). %
We then explore, mainly for the one variable case ($s=1$) and with
$r\in\mathbb{Z}$, some of the properties of the rational functions/polynomials
$Q_{n}^{r}(x)$, and list a few of them for low $r,n$.

We end the introduction with some comments on existing literature.
The finite product cases ($r\in\mathbb{N}_{0}$) of the identities
in Theorem \ref{thm:Main} appeared in \cite[Thm. 5.27]{FS}, where
they are seen as consequences of the Macdonald-Cheah formula for the
Poincaré/mixed Hodge polynomials of symmetric products of finite CW
complexes/algebraic varieties. This is recalled in Example \ref{exa:Macdonald-Cheah}.
On the other hand, the present approach is much more elementary, places
all $r\in\mathbb{R}$ in the same setting, and applies to several
variables.%

\subsection*{Acknowledgments}

I thank A. Nozad, J. Silva and A. Zamora for many conversations on
mixed Hodge structures and Serre polynomials of character varieties,
where some of the identities in this article were encountered and
distilled, and S. Mozgovoy for sparking my curiosity in plethysms.
Many thanks to Mark Wildon and Jimmy He for pointing out the relevance
of the cycle index in some statements, and to Peter Cameron, Yang-Hui
He and Richard Stanley for their interest, and for encouragement. 

\section{Plethystic Exponentials}

In this section, we recall the definitions and fundamental properties
of plethystic exponentials in multivariable formal power series rings,
with the goal of proving the fundamental product-sum identity in Equation
\eqref{eq:product-sum}. Usually, the plethystic exponential and logarithmic
functions are defined in the context of so-called $\lambda$-rings,
and related to Adams operations in algebra/algebraic topology. In
turn, these stem from the so-called plethysms, which have been widely
used in the theory of symmetric polynomials (see, for instance, \cite{LR,Kn,St}). 

As an application to Theoretical Physics, plethystic exponentials
are a fundamental ingredient of the famous DMVV formula for the orbifold
elliptic genus of symmetric products (\cite{DMVV}). More recently,
they were used for counting BPS gauge invariant operators in supersymmetric
gauge theories of $D$-branes probing Calabi-Yau singularities: the
\emph{plethystic program} of Feng, Hanany and He (\cite{BFHH,FHH}).
Here, these exponentials are introduced in elementary terms.

\subsection{One variable Plethystic exponentials}

We start by defining plethystic exponentials in a simple case: on
the ring $\mathbb{Q}\llbracket x\rrbracket$ of formal power series
in the variable $x$, with \emph{rational} \emph{coefficients}, written
as $f(x)=\sum_{n\geq0}a_{n}x^{n}$, $a_{n}\in\mathbb{Q}$. Recall
that $\mathbb{Q}\llbracket x\rrbracket$ has an obvious addition,
and the multiplication of such $f(x)$ with $g(x)=\sum_{n\geq0}b_{n}x^{n}$
is given by the Cauchy product: $f(x)\cdot g(x)=\sum_{n\geq0}(\sum_{k=0}^{n}a_{k}b_{n-k})x^{n}\in\mathbb{Q}\llbracket x\rrbracket.$
One can also compose formal power series $(f\circ g)(x)$ provided
$g(x)=\sum_{n\geq1}b_{n}x^{n}$, i.e., $g(x)$ has \emph{zero constant
term} ($b_{0}=0$):
\begin{equation}
(f\circ g)(x)=a_{0}+a_{1}(\sum_{n\geq1}b_{n}x^{n})+a_{2}(\sum_{n\geq1}b_{n}x^{n})^{2}+\cdots,\label{eq:composition}
\end{equation}
(this is well defined, as only a finite sum is involved in getting
the coefficient of a given $x^{n}$). Some important examples of formal
power series are the geometric and the (usual) exponential:
\begin{eqnarray*}
\frac{1}{1-x} & = & 1+x+x^{2}+\cdots+x^{n}+\cdots,\\
\exp(x) & = & 1+x+\frac{x^{2}}{2}+\cdots+\frac{x^{n}}{n!}+\cdots.
\end{eqnarray*}
For cases like these, which are actually convergent in some region,
we adopt the usual notations for analytic functions. Denote by $\mathbb{Q}^{0}\llbracket x\rrbracket:=x\mathbb{Q}\llbracket x\rrbracket$
the ideal consisting of series with zero constant term in $\mathbb{Q}\llbracket x\rrbracket$.
\begin{defn}
\label{def:PE-prod} Let $f(x)=\sum_{k\geq1}a_{k}x^{k}\in\mathbb{Q}^{0}\llbracket x\rrbracket$.
The \emph{harmonic operator} $\Psi:\mathbb{Q}^{0}\llbracket x\rrbracket\to\mathbb{Q}^{0}\llbracket x\rrbracket$,
is the linear map defined by: 
\[
\Psi[f](x)=\sum_{m\geq1}\frac{f(x^{m})}{m}\in\mathbb{Q}^{0}\llbracket x\rrbracket.
\]
The \emph{plethystic exponential} of $f$, denoted by $\P[f]$ is
the composition:
\begin{equation}
\P[f](x)=(\exp\circ\Psi[f])(x).\label{eq:PE-prod}
\end{equation}
\end{defn}

\begin{rem}
(1) Observe that we cannot define the harmonic operator on $\mathbb{Q}\llbracket x\rrbracket$,
since there is no meaning for the divergent \emph{harmonic series}:
$\sum_{m\geq1}\frac{1}{m}$. Also note that $\P[0]=1$.\\
(2) Often, for simplicity of notation, we write $\P(f(x))$ for $\P[f](x)$.
For example, if $f(x)=x+2x^{3}$, we write $\P[f](x)$ as $PE(x+2x^{3})$.\\
(3) Naturally, we can define $\P$ in $R\llbracket x\rrbracket$ for
any ring $R$ containing $\mathbb{Q}$. We will use \emph{real coefficients}
in subsection 2.2, and make occasional comments when using other coefficients. 
\end{rem}

\begin{prop}
\label{prop:propert}Let $f,g\in\mathbb{Q}^{0}\llbracket x\rrbracket$.
The plethystic exponential verifies:

(i) $\P[f](x)$ has constant term 1,

(ii) If $n\in\mathbb{N}$, $\P(x^{n})=\frac{1}{1-x^{n}}=1+x^{n}+x^{2n}+\cdots$,

(iii) $\P[f+g]=\P[f]\cdot\P[g]$,

(iv) $\P[-f]=(\P[f])^{-1}$.
\end{prop}

\begin{proof}
(i) Follows from Eq. \eqref{eq:composition}, as the constant term
of exp is 1. (ii) We evaluate the harmonic operator on $f(x)=x^{n}\in\mathbb{Q}^{0}\llbracket x\rrbracket$,
using the series expansion of the logarithm:
\[
\Psi[f](x)=\sum_{m\geq1}\frac{f(x^{m})}{m}=\sum_{m\geq1}\frac{x^{mn}}{m}=-\log(1-x^{n})=\log\left(\frac{1}{1-x^{n}}\right).
\]
Then, $\P[f](x)=\exp(\Psi[f](x))=1/(1-x^{n})$ as wanted. (iii) Follows
from the linearity of $\Psi$, and the usual properties of exp. (iv)
Immediate from (iii) and from $\P[0]=1$.
\end{proof}
\begin{example}
\label{exa:1}Using the rules above, some simple examples follow,
such as:

(i) With $n>0$, $\P(x^{n}-x^{2n})=\P(x^{n})/\P(x^{2n})=\frac{1-x^{2n}}{1-x^{n}}=1+x^{n}$,

(ii) $\P(\frac{x}{1-x})=\P(x+x^{2}+x^{3}+\cdots)=\prod_{k\geq1}\frac{1}{1-x^{k}}=1+x+2x^{2}+3x^{3}+5x^{4}+7x^{5}+\cdots$,

(iii) $\P\left(\frac{x(2-x)}{(1-x)^{2}}\right)=\prod_{k\geq1}\frac{1}{(1-x^{k})^{k+1}}=1+2x+6x^{2}+14x^{3}+33x^{4}+70x^{5}+149x^{6}+\cdots$.
\end{example}

The second and third examples above, which will show up later on,
relate to the famous Euler and Macmahon functions, respectively. They
are the generating functions for the sequence of partitions of natural
numbers $1,2,3,5,7,11,\cdots$, and the sequence $2,6,14,33,70,\cdots$
of partitions of $n\geq1$ objects colored using 2 colors (see §\ref{subsec:colored-partitions}
below).

These examples illustrate why the plethystic exponential plays such
an important role in the study of integer sequences, in the spirit
described in \cite{Ca}: indeed, when all the coefficients of $f$
are integers, $\P[f]$ corresponds precisely to the operator $S$
defined there, and is sometimes called the \emph{Euler transform}
(\cite{SP}). On the other hand, the proofs of Theorems \ref{thm:Main}
and \ref{thm:multivar} require the generalization to the multi-variable
plethystic exponential. 

\subsection{Multi-variable plethystic exponentials and partitions}

Consider the ring of formal power series with \emph{real coefficients}
in $s+1$ variables $x_{1},\cdots,x_{s}$ and $y$, $R:=\mathbb{R}\llbracket x_{1},\cdots,x_{s},y\rrbracket$,
and let $R^{0}\subset R$ be the ideal of series without constant
term. For example, when $s=1$, elements of $R^{0}$ can be written
alternatively as:
\[
f(x,y)=\sum_{j,k\geq0}\ a_{j,k}\ x^{j}\,y^{k}=\sum_{k\geq0}b_{k}(x)\,y^{k}.
\]
In the first expansion we have $a_{j,k}\in\mathbb{R}$ and $a_{0,0}=0$,
and in the second we have $b_{k}(x)\in\mathbb{R}\llbracket x\rrbracket$
and $b_{0}(x)\in\mathbb{R}^{0}\llbracket x\rrbracket$. When $s>1$,
to simplify the notation, we write $\underline{x}:=(x_{1},\cdots,x_{s})$,
$\mathbb{R}\llbracket\underline{x},y\rrbracket=\mathbb{R}\llbracket x_{1},\cdots,x_{s},y\rrbracket$,
and $\underline{x}{}^{m}$ means the $s$-tuple $(x_{1}^{m},\cdots,x_{s}^{m})$,
when $m\in\mathbb{N}$.
\begin{defn}
\label{def:PE-2var}Let $f(\underline{x},y)\in R^{0}$. The \emph{plethystic
exponential} of $f$, denoted by $\P[f]$ is the composition:
\[
\P[f]=\exp\circ\Psi[f]\in R,
\]
with the harmonic operator $\Psi:R^{0}\to R^{0}$, now given by $\Psi[f](\underline{x},y)=\sum_{m\geq1}\frac{f(\underline{x}^{m},y^{m})}{m}\,\in R^{0}.$
\end{defn}

We now need a computational tool for $\P[f]$, for functions $f$
of the form: 
\[
f(\underline{x},y)=g(\underline{x})\,y,\quad\quad g(\underline{x})\in\mathbb{R}\llbracket\underline{x}\rrbracket,
\]
which can be expressed in terms of partitions of natural numbers;
they are also certain \emph{cycle indices} for symmetric groups (see
Section 3.2). In \cite{FHH}, these are called ``$y$-inserted''
plethystic exponentials.

A partition of $n\in\mathbb{N}$, is given by the sum $n=n_{1}+\cdots+n_{m}$
where each $n_{j}$ is a positive integer. Alternatively, with $\mathbb{N}_{0}=\mathbb{N}\cup\{0\}$,
a partition $\lambda$ of $n$ is a map $\lambda:\{1,\cdots,n\}\mapsto\mathbb{N}_{0}$
such that $\lambda(j)$ is the number of parts of size $j$ (zero,
if there are no parts of size $j$). In this notation, we have $\sum_{j=1}^{n}j\lambda(j)=n$.
As usual, we write $\lambda\vdash n$ when $\lambda$ is a partition
of $n$.

We now show a fundamental \emph{product-sum formula} for plethystic
exponentials. %
{} Recall the multi-index notation: $\mathbf{k}=(k_{1},\cdots,k_{s})\in\mathbb{N}_{0}^{s}$,
and $\mathbf{x}^{\mathbf{k}}=x_{1}^{k}\cdots x_{s}^{k}$, for the
$s$ variables $\underline{x}=(x_{1},\cdots,x_{s})$.
\begin{prop}
\label{prop:Main2} Let $g(\underline{x})=\sum_{\mathbf{k}\in\mathbb{N}_{0}^{r}}a_{\mathbf{k}}\mathbf{x}^{\mathbf{k}}\in\mathbb{R}\llbracket\underline{x}\rrbracket$.
We have the equalities in $R=\mathbb{R}\llbracket\underline{x},y\rrbracket$
\begin{equation}
\P\left(g(\underline{x})\,y\right)=\prod_{\mathbf{k}\in\mathbb{N}_{0}^{s}}(1-y\,\mathbf{x}^{\mathbf{k}})^{-a_{\mathbf{k}}}=1+\sum_{n\geq1}\left(\sum_{\lambda\vdash n}\ \prod_{j=1}^{n}\frac{g(\underline{x}^{j})^{\lambda(j)}}{\lambda(j)!\,j^{\lambda(j)}}\right)y^{n},\label{eq:partition}
\end{equation}
where $\lambda(j)$ is the number of parts of $\lambda\vdash n$ of
each size $j\in\{1,2,\cdots,n\}$.
\end{prop}

\begin{proof}
Starting from the definition, we perform the following computation
(cf. also \cite{FNZ,FHH}):
\begin{eqnarray*}
\P[f](\underline{x},y)=\exp\left(\sum_{m\ge1}\frac{g(\underline{x}^{m})\,y^{m}}{m}\right) & = & \prod_{m\geq1}\exp\left(\frac{g(\underline{x}^{m})}{m}\,y^{m}\right)\\
 & = & \sum_{k_{1}\geq0}\sum_{k_{2}\geq0}\cdots\prod_{m\geq1}\,\frac{1}{k_{m}!}\left(\frac{g(\underline{x}^{m})}{m}\,y^{m}\right)^{k_{m}}.
\end{eqnarray*}
To get the coefficient of $y^{n}$, write $n=\sum_{k_{m}\geq0}mk_{m}$,
and we see that we are dealing with partitions $\lambda$ of $n$
with $k_{j}=\lambda(j)$ being the number of parts of size $j$ and
we get:\\[-3mm]
\[
\P[f](\underline{x},y)=1+\sum_{n\geq1}\left(\sum_{\lambda\vdash n}\ \prod_{j=1}^{n}\frac{g(\underline{x}^{j})^{\lambda(j)}}{\lambda(j)!\,j^{\lambda(j)}}\right)y^{n},
\]
which is the right hand side. To obtain the product form, note that
an analogous computation as in the proof of Proposition \ref{prop:propert}(ii)
gives $\P(y\,x_{1}^{k_{1}}\cdots x_{s}^{k_{s}})=\frac{1}{1-y\,\mathbf{x}^{\mathbf{k}}}$,
for every $k_{1},\cdots,k_{s},n\geq0$ not all zero. So, using Proposition
\ref{prop:propert}(iii) we transform the sum into a product :\\[-3mm]
\[
\P(\sum_{\mathbf{k}\in\mathbb{N}_{0}^{r}}a_{\mathbf{k}}\mathbf{x}^{\mathbf{k}}\,y)=\prod_{\mathbf{k}\in\mathbb{N}_{0}^{s}}(1-y\,\mathbf{x}^{\mathbf{k}})^{-a_{\mathbf{k}}}\in\mathbb{Z}\llbracket\underline{x},y\rrbracket,
\]
as wanted. Note that, if $-a_{\mathbf{k}}\in\mathbb{N}$, then $\P(a_{\mathbf{k}}\mathbf{x}^{\mathbf{k}}\,y)=(1-y\,\mathbf{x}^{\mathbf{k}})^{-a_{\mathbf{k}}}$
is actually a polynomial.
\end{proof}
\begin{rem}
\label{rem:Z-coeffs}(1) It is clear that the subset $1+\mathbb{R}^{0}\llbracket\underline{x}\rrbracket\subset\mathbb{R}\llbracket\underline{x}\rrbracket$
(series with 1 as constant term) is a commutative group with multiplication
of formal power series, and 1 as the identity. If we consider $(\mathbb{R}\llbracket\underline{x}\rrbracket,+,0)$
also as an abelian group one sees that $\P:\mathbb{R}^{0}\llbracket\underline{x}\rrbracket\to1+\mathbb{R}^{0}\llbracket\underline{x}\rrbracket$
is an \emph{isomorphism of abelian groups}. Its inverse, is called
naturally the \emph{plethystic logarithm} and it is also important
in enumerative combinatorics. For example, it allows the enumeration
of connected simple graphs, from the knowledge of the generating function
of all graphs (see \cite{Ha}).\\
(2) Even though Equations \eqref{eq:partition} appear to be make
sense only in $\mathbb{Q}\llbracket\underline{x},y\rrbracket$ or
$\mathbb{R}\llbracket\underline{x},y\rrbracket$, Proposition \ref{prop:Main2}
shows that if $f\in\mathbb{Z}^{0}\llbracket\underline{x},y\rrbracket\subset\mathbb{Q}^{0}\llbracket\underline{x},y\rrbracket$
(series with zero constant term), then $\P[f]$ is in $\mathbb{Z}\llbracket\underline{x},y\rrbracket$:
plethystic exponentials of series with \emph{integer coefficients},
still have integer coefficients.
\end{rem}

\section{Plethystic Exponentials of Rational Functions }

{} We come now to the proof of our main identities \eqref{eq:main-identity}
and \eqref{eq:multivar}. They follow from Theorem \ref{prop:Main2},
using the simple \emph{generalized} rational function:
\[
g(\underline{x})=(1-x_{1})^{r_{1}}\cdots(1-x_{s})^{r_{s}}\in\mathbb{R}\llbracket\underline{x}\rrbracket,
\]
for certain choices of \emph{real} numbers $r_{1},\cdots,r_{s}$.
We then explore some consequences of the identities, their relation
to colored partitions and the Euler and Macmahon functions, and to
the Macdonald/Cheah formula for symmetric products. We also observe
some of their consequences when used in conjunction with the so-called
$q$-binomial theorem and Heine's summation formula. We end with a
study of some simple properties of the polynomials $Q_{n}^{r}(x)$
in Equation \eqref{eq:char-pol-power}, which are evaluations of the
cycle index of $S_{n}$, and explicitly list some of them. 

\subsection{The main identities}

For $n\in\mathbb{N}$, let $S_{n}$ be the symmetric group on $n$
letters, that is, the group of permutations of the set $\{1,\cdots,n\}$.
This group acts linearly on a $n$ dimensional vector space, such
as $\mathbb{R}^{n}$, by permuting the elements of a fixed standard
basis. This way, we identify every element $\sigma\in S_{n}$ with
an element $M_{\sigma}\in GL(\mathbb{R}^{n})$ (the group of linear
automorphisms of $\mathbb{R}^{n}$ with identity $I_{n}$). Every
permutation $\sigma\in S_{n}$ can be written as a product $\sigma=\gamma_{1}\gamma_{2}\cdots\gamma_{m}$
where each $\gamma_{i}$ is a cycle, say of length $n_{i}$, and all
such cycles are disjoint. If we also include in this product all cycles
of length 1 (elements fixed by $\sigma$), then $n=n_{1}+n_{2}+\cdots$
is a partition of $n$. Conversely, any partition of $n$ determines
a unique conjugacy class of elements $\sigma\in S_{n}$.
\begin{thm}
\label{thm:multivariable} Fix $s\in\mathbb{N}$ and $r_{1},\cdots,r_{s}\in\mathbb{R}$.
Then, the following are equalities in $\mathbb{R}\llbracket\mathbf{x},y\rrbracket$:
\[
\P((1-x_{1})^{r_{1}}\cdots(1-x_{s})^{r_{s}}y)\ =\,1+\sum_{n\geq1}\sum_{\sigma\in S_{n}}\frac{y^{n}}{n!}\prod_{i=1}^{s}\det(I_{n}-x_{j}M_{\sigma})^{r_{i}}.
\]
\end{thm}

\begin{proof}
Given a permutation $\sigma\in S_{n}$, written as a product of disjoint
cycles $\sigma=\gamma_{1}\gamma_{2}\cdots\gamma_{m}$, the endomorphism
$I_{n}-xM_{\sigma}$ is a direct sum:
\[
I_{n}-xM_{\sigma}=(I_{n_{1}}-xM_{\gamma_{1}})\oplus\cdots\oplus(I_{n_{m}}-xM_{\gamma_{m}}).
\]
(with $M_{\gamma_{i}}$ considered as acting on a $n_{i}$-dimensional
subspace of $\mathbb{R}^{n}$). For a cycle $\gamma$ of length $k$,
we readily compute $\det(I_{k}-xM_{\gamma})=1-x^{k}$, so we see that
\[
\det(I_{n}-xM_{\sigma})^{r}=\prod_{j=1}^{n}(1-x^{j})^{r\sigma_{j}},
\]
where $\sigma_{j}$ is the number of cycles of $\sigma$ with length
$j\in\{1,\cdots,n\}$, and this determinant is invariant under conjugation
of $\sigma$ in $S_{n}$. Hence, in the multivariable case, and since
the number of $\sigma\in S_{n}$ with cycle type given by the partition
$\lambda\vdash n$, i.e. $(\lambda(1),\cdots,\lambda(n))=(\sigma_{1},\cdots,\sigma_{n})$,
is exactly $n!\prod_{j=1}^{n}\frac{1}{\sigma_{j}!\,j^{\sigma_{j}}}$,
we can sum over partitions to get:
\begin{eqnarray*}
1+\sum_{n\geq1}\frac{y^{n}}{n!}\sum_{\sigma\in S_{n}}\prod_{i=1}^{s}\det(I_{n}-x_{i}M_{\sigma})^{r_{i}} & = & 1+\sum_{n\geq1}\frac{y^{n}}{n!}\sum_{\sigma\in S_{n}}\prod_{j=1}^{n}\left[(1-x_{1}^{j})^{r_{1}}\cdots(1-x_{s}^{j})^{r_{s}}\right]^{\sigma_{j}}\\
 & = & 1+\sum_{n\geq1}\left(\sum_{\lambda\vdash n}\ \prod_{j=1}^{n}\frac{g(\underline{x}^{j})^{\lambda(j)}}{\lambda(j)!\,j^{\lambda(j)}}\right)y^{n}\\
 & = & \P(g(\underline{x})\,y),
\end{eqnarray*}
where $g(\underline{x})=(1-x_{1})^{r_{1}}\cdots(1-x_{s})^{r_{s}}$,
by Proposition \ref{prop:Main2}.
\end{proof}
We now complete the proof of theorems \ref{thm:Main} and \ref{thm:multivar}.
\begin{proof}
(of Theorem \ref{thm:Main}) Setting $x_{1}=\cdots=x_{s}$ and $r_{1}+\cdots+r_{s}=r$
in Theorem \ref{thm:multivariable} we are computing the plethystic
exponential of $f(x,y)=(1-x)^{r}y$. By the product formula of Proposition
\ref{prop:Main2}, this amounts to $\prod_{k\geq0}(1-yx^{k})^{-a_{k}}$
with $a_{k}$ being the coefficients of the series expansion of
\[
g(x)=(1-x)^{r}=\sum_{k\geq0}\binom{k-r-1}{k}x^{k},
\]
valid for every real $r$. When $r\in\mathbb{N}$, from $\binom{k-r-1}{k}=(-1)^{k}\binom{r}{k}$
we obtain Equation \eqref{eq:positive-r}, and the convergence of
the series follows from the estimate:
\[
\sum_{\sigma\in S_{n}}|\det(I_{n}-xM_{\sigma})|^{r}\leq\sum_{\sigma\in S_{n}}2^{nr}=n!2^{nr},
\]
valid when $|x|<1$. 
\end{proof}
For the proof of Equation \eqref{eq:multivar}, recall that we restrict
to integer coefficients.
\begin{proof}
(of Theorem \ref{thm:multivar}) From Theorem \ref{thm:multivariable},
the summation on the left hand side of Equation \eqref{eq:multivar}
equals $\P(\frac{(1-x_{1})\cdots(1-x_{r})}{(1-q_{1})\cdots(1-q_{s})}y)$.
Given the series expansion:
\[
\frac{(1-x_{1})\cdots(1-x_{r})}{(1-q_{1})\cdots(1-q_{s})}=\sum_{j_{1}=0}^{1}\cdots\sum_{k_{1}=0}^{\infty}\cdots(-x_{1})^{j_{1}}\cdots(-x_{r})^{j_{r}}q_{1}^{k_{1}}\cdots q_{s}^{k_{s}},
\]
the result follows immediately from the product formula, Proposition
\ref{prop:Main2}.
\end{proof}
\begin{rem}
We should point out that a special case of Theorem \ref{thm:multivariable}
had already been obtained in the study of symmetric products of algebraic
groups (see \cite[Lem. 4.2]{Si}); namely, for a sequence $r_{1},\cdots,r_{s}\in\mathbb{N}_{0}$:
\[
\P\left[y\,{\textstyle \prod_{j=1}^{s}}(1-x^{2j-1})^{r_{j}}\right]\ =\,1+\sum_{n\geq1}\sum_{\sigma\in S_{n}}\frac{y^{n}}{n!}\prod_{j=1}^{s}\det(I_{n}-x^{2j-1}M_{\sigma})^{r_{j}}.
\]
\end{rem}

\subsection{\label{subsec:colored-partitions}Applications and relations to other
identities}

We now present some relations of the main identities with other subjects,
and explore some consequences of these connections.
\begin{example}
(Generating function for colored partitions) We now prove Corollary
\ref{cor:residue}, which determines the residue of our basic generating
function
\[
\Phi^{r}(x,y):=1+\sum_{n\geq0}\sum_{\sigma\in S_{n}}\frac{y^{n}}{n!}\det(I_{n}-xM_{\sigma})^{r},\quad\quad r\in\mathbb{Z}
\]
at $y=1$. In fact, this is an immediate consequence of our main identity:
\[
\mathrm{Res}_{y=1}\Phi^{r}(x,y)=\mathrm{Res}_{y=1}\left((1-y)^{-1}\prod_{k\geq1}(1-yx^{k})^{-\binom{k-r-1}{k}}\right)=-\prod_{k\geq1}(1-x^{k})^{-\binom{k-r-1}{k}}.
\]
Two special cases of this formula are noteworthy, namely when $r=-1$
or $r=-2$ ($r=0$ is the geometric series $\Phi^{0}(x,y)=\frac{1}{1-y}$
with residue $-1$):
\begin{eqnarray*}
-\mathrm{Res}_{y=1}\Phi^{-1}(x,y) & = & \frac{1}{\phi(x)}=\sum_{n\geq0}p(n)x^{n}=1+x+2x^{2}+3x^{3}+5x^{4}+\cdots,\\
-\mathrm{Res}_{y=1}\Phi^{-2}(x,y) & = & \frac{M(x)}{\phi(x)}=1+2x+6x^{2}+14x^{3}+33x^{4}+70x^{5}+149x^{6}+\cdots,
\end{eqnarray*}
where $\phi(x):=\prod_{k\geq1}(1-x^{k})$ and $M(x)=\prod_{k\geq1}(1-x^{k})^{-k}$
are the famous Euler function and Macmahon function, respectively,
which are well known to be the generating functions for \emph{partitions}
and for \emph{plane partitions}, respectively. These are the sequences
A000041 and A000219 in \cite{OEIS}. Note also that $\frac{M(x)}{\phi(x)}=\prod_{k\geq1}(1-x^{k})^{-k-1}$
is actually the generating function for the number of \emph{colored
partitions of $n$ with 2 colors} (sequence A005380 in \cite{OEIS}).
More generally, a colored partition of $n\in\mathbb{N}$ with $c\in\mathbb{N}$
colors is defined as follows. Take a partition $n_{1}\geq n_{2}\geq\cdots\geq n_{m}\geq1$
of $n$ (so that $n=n_{1}+\cdots+n_{m}$) and consider each part $n_{j}$
as $n_{j}$ identical objects which are then colored with one of $c$
possible colors; we disregard the order of the colored objects. For
example, with $c=3=|\{b,w,r\}|$, for black, white and red colors,
colored partitions of $n=2$ are the 12 elements:
\[
bb,\quad bw,\quad br,\quad ww,\quad wr,\quad rr,\quad b+b,\quad b+w,\quad b+r,\quad w+w,\quad w+r,\quad r+r.
\]
Hence, Corollary \ref{cor:residue} follows from the well-known formula
for the generating function for colored partitions of $n$ with $c$
colors is (see \cite{OEIS}, sequence A217093, and \cite[exercise 7.99]{St2}):
\[
\prod_{k\geq1}(1-x^{k})^{-\binom{k-r-1}{k}}.
\]
\end{example}

\begin{example}
\label{exa:q-binomial}(the $q$-binomial theorem). In Theorem \ref{thm:multivar},
consider one variable $x$, and one $q$. We then get the infinite
product of the famous $q$-binomial theorem:
\begin{equation}
1+\sum_{n\geq0}\frac{y^{n}}{n!}\sum_{\sigma\in S_{n}}\frac{\det(I_{n}-xM_{\sigma})}{\det(I_{n}-qM_{\sigma})}=\prod_{k\in\mathbb{N}_{0}}\frac{1-xyq^{k}}{1-yq^{k}}=\sum_{n\geq0}\frac{(x;q)_{n}}{(q;q)_{n}}\,y^{n},\label{eq:q-binomial}
\end{equation}
where $(a;q)_{n}$ denote the (finite) Pochhammer symbols:
\[
(a;q)_{n}:=(1-a)(1-aq)\cdots(1-aq^{n-1}),\quad\quad(a;q)_{0}:=1.
\]
The right hand side of \eqref{eq:q-binomial} is also the \emph{basic
hypergeometric series} denoted by $_{1}\phi_{0}$. We conclude the
interesting formula:
\begin{equation}
\frac{1}{n!}\sum_{\sigma\in S_{n}}\frac{\det(I_{n}-xM_{\sigma})}{\det(I_{n}-qM_{\sigma})}=\frac{(x;q)_{n}}{(q;q)_{n}}=\prod_{k=1}^{n}\frac{1-xq^{k-1}}{1-q^{k}}.\label{eq:generalized-Molien}
\end{equation}
In turn, by letting $x=0$ in Equation \eqref{eq:generalized-Molien}
we get Molien's formula for the Hilbert series that captures the dimensions
of the graded components of the natural action of $S_{n}$ on the
polynomial ring in $n$ variables, $\mathbb{C}[t_{1},\cdots,t_{n}]$:
\[
\frac{1}{n!}\sum_{\sigma\in S_{n}}\frac{1}{\det(I_{n}-qM_{\sigma})}=\prod_{k=1}^{n}\frac{1}{1-q^{k}}.
\]
We thank Jimmy He for suggesting to look for analogies between our
identities and other $q$-series identities. 
\end{example}

\begin{example}
(Heine's summation) Consider now two variables $x$, and one $q$.
We then obtain:
\[
1+\sum_{n\geq0}\frac{y^{n}}{n!}\sum_{\sigma\in S_{n}}\frac{\det[(I_{n}-x_{1}M_{\sigma})(I-x_{2}M_{\sigma})]}{\det(I_{n}-qM_{\sigma})}=\P\left(\frac{(1-x_{1})(1-x_{2})}{1-q}y\right)
\]
which is precisely the infinite product in the so-called Heine's summation
theorem, in the form:
\[
\sum_{n\geq0}\prod_{k=1}^{n}\frac{(1-x_{1}^{-1}q^{k-1})(1-x_{2}^{-1}q^{k-1})}{(1-q^{k})(1-yq^{k-1})}\,(x_{1}x_{2}y)^{n}=\prod_{k\in\mathbb{N}_{0}}\frac{(1-yx_{1}q^{k})(1-yx_{2}q^{k})}{(1-yq^{k})(1-yx_{1}x_{2}q^{k})}
\]
since this is the basic hypergeometric series $_{2}\phi_{1}(\frac{1}{x_{1}},\frac{1}{x_{2}};\,y;\,q;\,yx_{1}x_{2})$
(see, e.g., \cite[(1.5.1)]{GR}). The generalization of these identities
to other cases, and their interesting relation with more general hypergeometric
series will be dealt in forthcoming work (\cite{F}). 
\end{example}

\begin{example}
\label{exa:Macdonald-Cheah}(Relation with Macdonald/Cheah formula
for symmetric products). Let $X$ be a $d$-dimensional finite CW
complex, with Poincaré polynomial $P_{X}(x)=\sum_{k=0}^{d}b_{k}(X)\,x^{k}$,
where $b_{k}(X)$ are its Betti numbers. The $n$th symmetric product
$\sym^{n}X:=X^{n}/S_{n}$ is the quotient of the $n$-fold cartesian
product under the action of $S_{n}$ by permutation of the variables.
The famous Macdonald formula reads (see \cite[pg. 563]{Mac}):
\begin{equation}
1+\sum_{n\geq1}P_{\sym^{n}X}(-x)\,y^{n}=\prod_{k=0}^{d}\left(1-x^{k}y\right)^{(-1)^{k+1}b_{k}(X)},\label{eq:Macdonald}
\end{equation}
and we see that the right hand side is naturally a plethystic exponential:
from Proposition \ref{prop:Main2}, it is precisely $\P(P_{X}(-x)y)$.
This formula was generalized for the mixed Hodge polynomials $\mu_{X}$
of symmetric products of complex algebraic varieties by J. Cheah in
\cite{Ch}, which reads, in our notation:
\[
1+\sum_{n\geq1}\mu_{\sym^{n}X}(t,u,v)\,y^{n}=\P(\mu_{X}(-t,u,v)\,y).
\]
More recently, by studying the space of commuting elements in a compact
Lie group $G$ (see \cite{Sf}), it was shown that:
\begin{equation}
P_{T^{r}/W}(x)=\frac{1}{|W|}\sum_{\sigma\in W}\det(I+x\sigma)^{r},\label{eq:charvar}
\end{equation}
where $W$ is the Weyl group of $G$ acting by reflections on the
dual of the Lie algebra of a maximal torus $T$, and $I$ is the identity
automorphism. An analogous study for mixed Hodge structures on character
varieties of abelian groups led to a proof of the positive integer
case $(r\geq1)$ of Theorem \ref{thm:Main} in \cite{FS}, by combining
equations \eqref{eq:Macdonald} and \eqref{eq:charvar}, in the case
of the unitary group $G=U(n)$, with maximal torus $T=(S^{1})^{n}$
and Weyl group $S_{n}$ (and using the deformation retraction from
these varieties to the quotient $T^{r}/W$, see \cite{FL}). 
\end{example}

\subsection{The cycle index and average powers of characteristic polynomials
of permutations}

A crucial ingredient in the proof of Theorem \ref{thm:multivariable}
was the cycle type of a permutation; so, there is a natural connection
between plethystic exponentials of rational functions and cycle indices
of $S_{n}$. Given $n$ variables $p_{1},\cdots,p_{n}$, the cycle
index of $S_{n}$ is defined as the degree $n$ polynomial, with $\mathbb{Q}$
coefficients: 
\[
Z_{n}(p_{1},\cdots,p_{n})=\frac{1}{n!}\sum_{\sigma\in S_{n}}\ \prod_{j=1}^{n}p_{j}^{\sigma_{j}}\quad\in\mathbb{Q}[p_{1},\cdots,p_{n}],
\]
where $\sigma_{j}\geq0$ is the number of cycles of $\sigma\in S_{n}$
of each size $j\in\{1,2,\cdots,n\}$. Now, if $g(\underline{x})\in\mathbb{R}\llbracket\underline{x}\rrbracket$
is a formal power series in $\underline{x}=(x_{1},\cdots,x_{n})$,
we denote by $Z_{n}[g(\underline{x})]$ the substitution 
\[
Z_{n}(g(\underline{x}),g(\underline{x}^{2}),\cdots,g(\underline{x}^{n})),
\]
and call it the \emph{evaluation of $Z_{n}$ at $g(\underline{x})$};
this is generally \emph{not a polynomial} in $x_{1},\cdots,x_{s}$,
but another formal power series in $\mathbb{R}\llbracket\underline{x}\rrbracket$.
We thank M. Wildon for noticing the relevance of the evaluation of
$Z_{n}$ in the context of the more general Proposition \ref{prop:Main2}
(\cite{Wi}) (leading also to the treatment of other results in this
section). In fact, since the number of permutations with the cycle
type of the partition $\lambda\vdash n$ equals $n!\prod_{j=1}^{n}\frac{1}{\lambda(j)!\,j^{\lambda(j)}}$,
the evaluation of $Z_{n}$ on an arbitrary $g(\underline{x})\in\mathbb{R}\llbracket\underline{x}\rrbracket$,
gives the following, in terms of partitions:
\[
Z_{n}[g(\underline{x})]=\frac{1}{n!}\sum_{\sigma\in S_{n}}g(\underline{x})^{\sigma_{1}}g(\underline{x}^{2})^{\sigma_{2}}\cdots g(\underline{x}^{n})^{\sigma_{n}}=\sum_{\lambda\vdash n}\ \prod_{j=1}^{n}\frac{g(\underline{x}^{j})^{\lambda(j)}}{\lambda(j)!\,j^{\lambda(j)}}\quad\in\mathbb{R}\llbracket\underline{x}\rrbracket.
\]
So, passing from partitions to permutations in Proposition \ref{prop:Main2}
gives the elegant relation:
\begin{equation}
\P(g(\underline{x})y)=1+\sum_{n\geq1}Z_{n}[g(\underline{x})]\,y^{n},\label{eq:cycle-index}
\end{equation}
for \emph{any} formal power series $g(\underline{x})\in\mathbb{R}\llbracket\underline{x}\rrbracket$
in an arbitrary number of variables $\underline{x}=(x_{1},\cdots,x_{s})$.
In the language of the \emph{Plethystic Program }of \cite{FHH}, this
says that ``$y$-inserted'' plethystic exponentials are just the
generating functions of the evaluations of cycle indices of $S_{n}$. 

Now fix a vector $\mathbf{r}=(r_{1},\cdots,r_{s})\in\mathbb{R}^{s}$
and consider the formal power series:
\[
Q_{n}^{\mathbf{r}}(\underline{x}):=\frac{1}{n!}\sum_{\sigma\in S_{n}}\det(I_{n}-x_{1}M_{\sigma})^{r_{1}}\cdots\det(I_{n}-x_{s}M_{\sigma})^{r_{s}}=Z_{n}[(1-x_{1})^{r_{1}}\cdots(1-x_{s})^{r_{s}}].
\]
 To prove Corollary \ref{cor:recursion}, we just need to use the
analogous property of the cycle index of the symmetric group $Z_{n}$.
For the benefit of the reader, we reproduce the argument, that boils
down to taking the derivative of Equation \eqref{eq:cycle-index}
with respect to $y$:
\begin{eqnarray*}
\sum_{n\geq1}n\,Z_{n}[g(\underline{x})]\,y^{n-1} & = & \frac{\partial}{\partial y}\left(\exp\left[\sum_{k\geq1}g(\underline{x}^{k})\frac{y^{k}}{k}\right]\right)\\
 & = & \left(1+\sum_{m\geq1}Z_{m}[g(\underline{x})]\,y^{m}\right)\left(\sum_{k\geq1}g(\underline{x}^{k})y^{k-1}\right),
\end{eqnarray*}
which implies, upon picking the $y^{n-1}$ term:
\[
Z_{n}[g(\underline{x})]=\frac{1}{n}\sum_{k=1}^{n}g(\underline{x}^{k})\,Z_{n-k}[g(\underline{x})].
\]
In particular, 
\begin{equation}
Q_{n}^{\mathbf{r}}(\underline{x})=\frac{1}{n}\sum_{k=1}^{n}(1-x_{1}^{k})^{r_{1}}\cdots(1-x_{s}^{k})^{r_{s}}\,Q_{n-k}^{\mathbf{r}}(\underline{x}),\label{eq:recursion}
\end{equation}
which shows Corollary \ref{cor:recursion}. Consider now the one variable
case, $s=1$, and let
\[
Q_{n}^{r}(x):=\frac{1}{n!}\sum_{\sigma\in S_{n}}\det(I_{n}-xM_{\sigma})^{r},
\]
be the average $r$th power of permutation polynomials for a given
matrix size $n\in\mathbb{N}$. Then, the recurrence relation \eqref{eq:recursion}
becomes $Q_{n}^{r}(x)=\frac{1}{n}\sum_{k=1}^{n}(1-x^{k})^{r}Q_{n-k}^{r}(x),$
for any $r\in\mathbb{R}$, from which we obtain the following explicit
expressions. 
\begin{prop}
For $n=1,2,3$ and $4$, and all $r\in\mathbb{R}$, we have:\\
(i) $Q_{1}^{r}(x)=(1-x)^{r}$\\
(ii) $Q_{2}^{r}(x)=\frac{1}{2}\left((1-x)^{2r}+(1-x^{2})^{r}\right)$\\
(iii) $Q_{3}^{r}(x)=\frac{1}{6}(1-x)^{3r}+\frac{1}{2}(1-x^{2})^{r}(1-x)^{r}+\frac{1}{3}(1-x^{3})^{r}$\\
(iv) $Q_{4}^{r}(x)=\frac{1}{24}(1-x)^{4r}+\frac{1}{4}(1-x^{4})^{r}+\frac{1}{8}(1-x^{2})^{2r}+\frac{1}{3}(1-x)^{r}(1-x^{3})^{3r}+\frac{1}{4}(1-x^{2})^{r}(1-x)^{2r}$.
\end{prop}

In the positive integer case $r\in\mathbb{N}$, the functions $Q_{n}^{r}(x)$
are one variable \emph{polynomials}, and we derive some of their general
properties as follows.
\begin{prop}
For $r,n\geq1$, $Q_{n}^{r}(0)=1$ and the degree of $Q_{n}^{r}(x)$
is $\leq nr$. Moreover: \\
(i) $Q_{n}^{r}(x)$ is divisible by $Q_{1}^{r}(x)=(1-x)^{r}$, in
particular $Q_{n}^{r}(1)=0$.\\
(ii) $Q_{n}^{r}(x)$ has integer coefficients.\\
(iii) When $r$ is even, $Q_{n}^{r}(x)$ is \emph{palindromic} (hence
monic of degree $nr$).
\end{prop}

\begin{proof}
Since the degree of the permutation polynomial $\det(I_{n}-xM_{\sigma})$
is $n$, and $|S_{n}|=n!$, the first sentence is clear. (i) For any
permutation $\sigma\in S_{n}$, the matrix $M_{\sigma}\in GL(n,\mathbb{C})$
always has an eigenvalue equal to 1; hence, for $r\geq1$, $Q_{1}^{r}(x)=(1-x)^{r}$
divides $\det(I_{n}-xM_{\sigma})^{r}$. (ii) This is an immediate
consequence of the product expansion of plethystic exponentials and
Remark \ref{rem:Z-coeffs}(2). (iii) A polynomial $p(x)$ is palindromic
if and only if 
\[
p(x)=x^{\deg(p)}p({\textstyle \frac{1}{x}}).
\]
Since every cycle $\delta\in S_{n}$ of length $k$ satisfies $\det(I_{k}-xM_{\delta})^{r}=(1-x^{k})^{r}$,
which is a palindromic polynomial when $r\in2\mathbb{N}$, and every
permutation is a product of cycles, the result follows from the fact
that the product of palindromic polynomials is palindromic.
\end{proof}
{} Finally, for low values of $r\in\mathbb{N}$ and arbitrary $n$,
the following is obtained by simple computations.
\begin{prop}
For $r=1,2$ and $3$, we have: \\
(i) $Q_{n}^{1}(x)=(1-x)$, for all $n\geq1$.\\
(ii) $Q_{n}^{2}(x)=(1-x)^{2}\cdot\sum_{k=0}^{n-1}x^{2k}$, for all
$n\geq1$.\\
(iii) $Q_{2}^{3}(x)=(1-x)^{3}(1+3x^{2})$\\
(iv) $Q_{n}^{3}(x)=(1-x)^{3}\left(\binom{n+1}{2}x^{2n-2}+\binom{n}{2}x^{2n-4}+\sum_{k=0}^{n-3}(k+1)x^{2k}\right)$,
for $n\geq3$.
\end{prop}

\end{document}